\documentclass[a4paper,12pt]{amsart}
\usepackage{amsmath,amssymb,latexsym,amsfonts,amscd}

\title{On the geography of Gorenstein minimal 3-folds of general type}
\author{Meng Chen and Christopher D. Hacon}
\address{\rm Institute of Mathematics, School of Mathematical Sciences, Fudan University,
Shanghai, 200433, People's Republic of China}
\email{mchen@fudan.edu.cn}

\address{\rm Department of Mathematics, University of Utah, 155 South 1400 East, Room 233, Salt Lake City, Utah 84112-0090, USA}
\email{hacon@math.utah.edu}

\thanks{The first author was supported both by Program for New Century Excellent Talents in
University (Ministry of Education, PR China) and by the National
Natural Science Foundation of China (no: 10131010). The second
author was partially supported by NSF research grant
  no: 0456363}
\newcommand{\ot}{{\otimes}}
\newcommand{\OO}{{\mathcal{O}}}

\newcommand{\roundup}[1]{\ulcorner{#1}\urcorner}

\newcommand{\Pic}{{\operatorname{Pic}}}
\newcommand{\Z}{{\mathbb {Z}}}

\newtheorem{thm}{Theorem}[section]
\newtheorem{lem}[thm]{Lemma}

\newtheorem{prop}[thm]{Proposition}

\theoremstyle{definition}

\newtheorem{setup}[thm]{}
\newtheorem{question}[thm]{Question}
\newtheorem{exmp}[thm]{Example}

\newtheorem{rem}[thm]{Remark}
\theoremstyle{remark}

\begin{document}
\begin{abstract}
Let $X$ be a minimal projective Gorenstein 3-fold of general type.
We give two applications of an inequality between $\chi (\omega _X)$ and $p_g(X)$:

1) Assume that the canonical map $\Phi_{|K_X|}$ is of fiber type. Let
$F$ be a smooth model of a generic irreducible component in the
general fiber of $\Phi_{|K_X|}$. Then the
birational invariants of $F$ are bounded from above.

2) If $X$ is nonsingular, then $c_1^3\leq \frac 1 {27}
c_1c_2+\frac{10}{3}$.
\end{abstract}
\maketitle
\pagestyle{myheadings} \markboth{\hfill Meng Chen and Christopher
D. Hacon\hfill}{\hfill On the geography of Gorenstein minimal
3-folds\hfill}
\section{\bf Introduction}

We work over the complex number field ${\mathbb C}$.

The main purpose of this note is to study the geometry of Gorenstein minimal $3$-folds $X$ of general type. We improve the inequality $\chi (\omega _X)\leq 2p_g(X)$ (see Proposition \ref{right} for a precise statement), and we show how this
leads to several applications which we explain below:

First, we improve the main theorem in
\cite{chen2005}:

\begin{thm}\label{main} Let $X$ be a minimal projective Gorenstein 3-fold of
general type. Assume that the canonical map $\Phi_{|K_X|}$ is of
fiber type. Let $F$ be a smooth model of a generic irreducible
component in the general fiber of $\Phi_{|K_X|}$. Then the
invariants of $F$ are bounded from above as follows:

(1) if $F$ is a curve, then $g(F)\le 487$;

(2) if $F$ is a surface, then $p_g(F)\le 434$.
\end{thm}

\begin{rem} 1) Theorem \ref{main} was verified by the first author
in \cite{chen2005} under the assumption that $p_g(X)$ is
sufficiently large.

2) When $\Phi_{|K_X|}$ is generically finite, the generic degree
is bounded from above by the second author in \cite{hacon2004}.

3) In the surface case, the corresponding boundedness theorem was
proved by Beauville in \cite{beauville}.

4) The numerical bounds in the above theorem might be far from
sharp.
\end{rem}

Our second application is an inequality of Noether type between
$c_1$ and $c_2$ which improves the main theorem
of \cite{shin}.

\begin{thm}\label{inequality1}
Let $X$ be a nonsingular projective minimal 3-fold of general
type. Then the following inequality holds:
$$K_X^3\ge \frac{8}{9}\chi(\omega_X)-\frac{10}{3}, \text{or equivalently}$$
$$c_1^3\leq \frac{1}{27}c_1c_2+\frac{10}{3}.$$
\end{thm}

\bigskip

Chen is grateful to De-Qi Zhang for pointing out an inequality (see
the proof of Lemma 2.1(3) in \cite{ZDQ}) similar to the one in
Proposition \ref{right} and for an effective discussion.

\section{\bf Proof of Theorem \ref{main}}

Throughout this note, {\it a minimal 3-fold} $X$ is one with nef canonical
divisor $K_X$ and with only ${\mathbb Q}$-factorial terminal
singularities.

\begin{setup}{\bf Notations and the set up.}\label{notation}
Let $X$ be a minimal projective 3-fold of general type. Since we
are discussing the behavior of the canonical map, we may assume
$p_g(X)\ge 2$. Denote by $\varphi_1$ the canonical map which is
usually a rational map. Take the birational modification $\pi:
X'\longrightarrow X$, which exists by Hironaka's big theorem, such
that

(i) $X'$ is smooth;

(ii) the movable part of $|K_{X'}|$ is base point free;

(iii) there exists a canonical divisor $K_X$ such that
$\pi^*(K_X)$ has support with only normal
crossings.

Denote by $h$ the composition $\varphi_1\circ\pi$. So $h:
X'\longrightarrow W'\subseteq{\mathbb P}^{p_g(X)-1}$ is a morphism.
Let $h: X'\overset{f}\longrightarrow B\overset{s}\longrightarrow
W'$ be the Stein factorization of $h$. We can write
$$K_{X'}=\pi^*(K_X)+E=S+Z,$$
where $S$ is the movable part of $|K_{X'}|$, $Z$ is the fixed part
and $E$ is an effective ${\mathbb Q}$-divisor which is a sum of
distinct exceptional divisors.

If $\dim\varphi_1(X)<3$, $f$ is a called an {\it induced fibration
of} $\varphi_1$. If $\dim\varphi_1(X)=2$, a general fiber $F$ of
$f$ is a smooth curve $C$ of genus $g:=g(C)\ge 2$. If
$\dim\varphi_1(X)=1$, a general fiber $F$ of $f$ is a smooth
projective surface of general type. Denote by $F_0$ the smooth minimal
model of $F$ and by $\sigma: F\longrightarrow F_0$ the smooth blow down
map. Denote by $b$ the genus of the base curve $B$.
\end{setup}

\begin{prop}\label{right} Let $V$ be a smooth projective 3-fold of general
type with $p_g(V)>0$. Then $\chi(\omega_V)\le p_g(V)$ unless a generic
irreducible component in the general fiber of the Albanese morphism
is a surface
$V_y$ with $q(V_y)=0$, in which case one has the inequality
$$ \chi(\omega_V)\le(1+\frac 1 {p_g(V_y)})p_g(V).$$
\end{prop}
\begin{proof} Since $\chi(\omega_V)=p_g(V)+q(V)-h^2(\OO _V)-1$,
the result is clear if $q(V)\leq 1$. So assume that $q(V)\geq 2$.
Let $a:V\to Y$ be the Stein
factorization of the Albanese morphism $V\to A(V)$. By the proof of
Theorem 1.1 in \cite{hacon2004}, one sees that we may assume that
${\rm dim} Y=1$ and hence $Y$ is a smooth curve.
Recall also that
by \cite{hacon2004}, $p_g(V)\geq \chi (a_* \omega _V)$.
Let $y\in Y$ be a general point and $V_y$ the corresponding fiber.
$V_y$ is a smooth surface of general type. If $q(V_y)>0$, then proceeding as
in \cite{hacon2004}, one sees that $\chi (R^1a_* \omega _V)=
\chi (R^1a_* \omega _{V/Y}\ot \omega _Y)$. Since the genus of $Y$ is $q(V)$,
and ${\rm deg} R^1a_* \omega _{V/Y}\geq 0$, one sees by an easy Riemann-Roch
computation that
$$\chi (R^1a_* \omega _V)\geq (q(V)-1)q(V_y).$$
Recall that $R^2a_* \omega _V \cong \omega _Y$ and so
$$\chi (\omega _V)=\chi (a_* \omega _V)-\chi (R^1a_* \omega _V)+\chi (R^2a_* \omega _V)
\leq \chi (a_* \omega _V)\leq p_g(V)$$ whenever $q(V_y)>0$.

We may therefore assume that $q(V_y)=0$. Notice that by \cite{Kollar},
the sheaf $R^1a_* \omega _V$ is torsion free. Since its rank is given by $h^1(\omega _{V_y})=q(V_y)=0$, we have that $R^1a_* \omega _V=0$.
Therefore, by a similar Riemann-Roch
computation, one sees that $\chi (a_* \omega _V)\geq
(q(V)-1)p_g(V_y)$ and so $$\chi (\omega _V)=\chi (a_* \omega
_V)+q(V)-1\leq \chi (a_* \omega _V)(1+\frac 1 {p_g(V_y)}) \leq
p_g(V)(1+\frac {1}{p_g(V_y)}).$$
\end{proof}
\begin{exmp} Let $S$ be a minimal surface of general type admitting a
$\mathbb{Z} _2$ action
such that $q(S)=0$, $p_g(S)=1$ and $p_g(S/\mathbb{Z} _2)=0$ (cf.
(2.6) of \cite{hacon2003}). Let $C$ be a curve admitting
a fixed point free $\mathbb{Z} _2$ action and let $B=C/\mathbb{Z} _2$.
Assume that the genus of $B$ is $b\geq 2$.
Let $V=S\times C/\mathbb{Z} _2$ be the quotient by the induced
diagonal action. Then $V$ is minimal, Gorenstein of general type
such that $p_g(V)=b-1$, $q(V)=b$ and $h^2(\OO _V)=0$.
It follows that $\chi (\omega _V)=2b-2=(1+1/p_g(V_y))p_g(V)$.


This example shows that the above proposition is close to being optimal.
\end{exmp}
\begin{lem}\label{pseudo} Let $X$ be a minimal 3-fold of general type. Suppose
$\dim\varphi_1(X)=1$. Keep the same notations as in \ref{notation}.
Replace $\pi:X'\longrightarrow X$, if necessary, by a further
birational modification (we still denote it by $\pi$). Then
$$\pi^*(K_X)|_F-\frac{p_g(X)-1}{p_g(X)}\sigma^*(K_{F_0})$$
is pseudo-effective.
\end{lem}
\begin{proof} One has an induced fibration $f:X'\longrightarrow
B$.

{\bf Case 1}. If $b>0$, we may replace $\pi$ by a new one as in the
proof of Lemma 2.2 of \cite{chenMA} such that $\pi^*(K_X)|_F\sim
\sigma^*(K_{F_0})$. In fact, since the fibers of $\pi$ are
rationally connected\footnote{Shokurov (\cite{Sho}) proved that if
the pair $(X, \Delta)$ is klt and the MMP holds, then the fibres of
the exceptional locus are always rationally chain connected.
Furthermore, the second author and M$^{\rm c}$Kernan (see
\cite{H-M}) have recently extended Shokurov's result to any
dimension and without assuming MMP.} and $b>0$, it follows that
$f:X'\to B$ factors through a morphism $f_1:X\to B$. But since $X$
is minimal and terminal, it follows that a general fiber $X_b$ of
$f_1$ is a smooth minimal surface of general type and hence it can
be identified with $F_0$. It is now clear that $\pi^*(K_X)|_F\sim
\sigma^*(K_{F_0})$.

Thus it
suffices to consider the case $b=0$.

{\bf Case 2}. If $p_g(X)=2$, the lemma was verified in section 4 (at
page 526 and page 527) in \cite{chen2003}. If $p_g(X)\ge 3$, one may
refer to Lemma 3.4 in \cite{chenMA}.
\end{proof}

\begin{prop}\label{inequality} Let $X$ be a Gorenstein minimal
projective 3-fold of general type. Let $d:=\dim \varphi _1 (X)$.
The following inequalities
hold:

(1) If $d=2$, then $K_X^3\ge
\roundup{\frac{2}{3}(g(C)-1)}(p_g(X)-2).$

(2) If $d=1$, then $K_X^3\ge
(\frac{p_g(X)-1}{p_g(X)})^2K_{F_0}^2(p_g(X)-1).$
\end{prop}
\begin{proof}
The inequality (1) is due to Theorem 4.1(ii) in \cite{chenJMSJ}.

Suppose now that $d=1$. We may write
$$\pi^*(K_X)\sim S+E_{\pi}$$
where $S\equiv tF$ with $t\ge p_g(X)-1$ and $E_{\pi}$ is an
effective divisor.

Thus we have

\begin{eqnarray*}
K_X^3&=&\pi^*(K_X)^3\ge (\pi^*(K_X)^2\cdot F)(p_g(X)-1)\\
&\ge& (\frac{p_g(X)-1}{p_g(X)})^2\sigma^*(K_{F_0})^2(p_g(X)-1)
\end{eqnarray*}
where Lemma \ref{pseudo} has been applied to derive the second
inequality above.
\end{proof}

\begin{setup}\label{p1}{\bf Proof of Theorem \ref{main}.}
The Miyaoka-Yau inequality (cf. \cite{miyaoka}) says
$$K_X^3\le 72\chi(\omega_X).$$

\begin{quote}
(**)\ \  Denote by $V$ a smooth model of $X$. Assume that a
generic irreducible component in the general fiber of the Albanese
morphism is a surface $V_y$ with $q(V_y)=0$ and $p_g(V_y)=1$.
Because $p_g(V)=p_g(X)\ge 2$, we see that the canonical map of $V$
maps $V_y$ to a point. This means $\dim\varphi_1(X)=1$, i.e.
$|K_X|$ is composed with a pencil.
Thus, one sees that in this special situation, the Stein factorization of the Albanese map is the fibration induced by $\varphi_1$.
So $p_g(F)=p_g(V_y)=1$.
\end{quote}

(1) Assume $\dim\varphi_1(X)=2$. The above argument implies that
$\chi(\omega_X)\leq \frac{3}{2}p_g(X)$ and so by proposition
\ref{inequality} and an easy computation, one sees that $g(C)\le
487$. Furthermore $g(C)\le 109$ whenever $p_g(X)$ is sufficiently
big.

(2) Assume $\dim\varphi_1(X)=1$. When $b>0$, we have $p_g(F)\le
38$ by (both 1.4 and Theorem 1.3 in) \cite{chen2005}. So we only
need to study the case $b=0$.

Suppose $p_g(F)\ge 2$. Then one has $\chi(\omega_X)\le
\frac{3}{2}p_g(X)$ by argument (**) and Proposition \ref{right}.
The Miyaoka-Yau inequality yields $K_X^3\le 72\chi(\omega_X)\le
108 p_g(X).$ Again by Propositions \ref{right} and
\ref{inequality}, we have
$$K_{F_0}^2\le 108 (\frac{p_g(X)}{p_g(X)-1})^3\le 864.$$
Also $K_{F_0}^2\le 108$ whenever $p_g(X)$ is sufficiently big. Taking into
account the Noether inequality $K_{F_0}^2\ge 2p_g(F)-4$, we get
$p_g(F)\le 434$. This concludes the proof.
\end{setup}

\section{\bf A Noether type inequality between $c_1$ and $c_2$}
\begin{setup}\label{known}{\bf A known inequality.} Let $X$ be a nonsingular
projective minimal 3-fold of general type. We have a sharp
inequality
$$K_X^3\ge \frac{4}{3}p_g(X)-\frac{10}{3}$$
which was first proved in \cite{chenmrl} under the assumption
$K_X$ being ample. The general case was recently proved in \cite{CCZ}.
\end{setup}

\begin{setup}{\bf Proof of Theorem \ref{inequality1}}\end{setup}
\begin{proof} Note that since $K_X^3>0$ is an even integer, the Theorem clearly holds for $\chi (\omega _X)\leq 6$. Therefore, we may assume that $\chi (\omega _X)>0$.

{\bf Case 1.} $p_g(X)>0$.

According to Proposition \ref{right}, we have an inequality
$\chi(\omega_X)\leq \frac{3}{2}p_g(X)$ unless a generic
irreducible component in the general fiber of the Albanese
morphism is a surface $V_y$ with $q(V_y)=0$ and $p_g(V_y)=1$.

So in the general case by \ref{known} one has the inequality
$$K_X^3\ge \frac{8}{9}\chi(\omega_X)-\frac{10}{3}$$
or equivalently,
$$c_1^3\le \frac{1}{27}c_1c_2+\frac{10}{3}.$$

In the exceptional case with $p_g(X)>1$, the argument (**) in
\ref{p1} says that $|K_X|$ is composed with a pencil of surfaces
and $\varphi_1$ generically factors through the Albanese map. Thus
$X$ is canonically fibred by surfaces with $q(V_y)=0$ and
$p_g(V_y)=1$. According to Theorem 4.1(iii) in \cite{chenJMSJ},
one has $K_X^3\ge 2p_g(X)-4$. Since by Proposition \ref{right}
$\chi(\omega_X)\le 2p_g(X)$, one has the stronger
inequality $K_X^3\ge \chi(\omega_X)-4$.

In the exceptional case with $p_g(X)=1$, by Proposition
\ref{right}, one has $\chi(\omega_X)\le 2$ and so the inequality in
Theorem \ref{inequality} holds.

{\bf Case 2.}  $p_g(X)=0$.

We can not rely on \ref{known} in this case. Since
$\chi(\omega_X)>0$, one has $q(X)>1$. Thus we can study the
Albanese map. Let $a:X\longrightarrow Y$ be the Stein
factorization of the Albanese morphism. We claim that $\dim(Y)=1$.
In fact, if $\dim (Y)\ge 2$, then the Proof of Theorem 1.1 of
\cite{hacon2004} shows $p_g(X)\ge \chi(a_* \omega_X)\ge \chi(\omega_X)>0$, a
contradiction.

So we have a fibration $a:X\longrightarrow Y$ onto a smooth curve
$Y$ with $g(Y)=q(X)>1$. Denote by $F$ a general fiber of $a$. If
$p_g(F)>0$, then the Proof of Theorem 1.1 of \cite{hacon2004} also
shows $0=2p_g(X)\ge \chi(\omega_X)>0$, which is also a contradiction. Thus
one must have $p_g(F)=0$. Because $F$ is of general type, one has
$q(F)=0$. Therefore, the sheaves $a_*\omega _X $ and $R^1a_*\omega _X $
have rank $h^0(\omega _F)=p_g(F)=0$ and $h^1(\omega _F)=q(F)=0$.
Since, by \cite{Kollar}, they are torsion free, it follows that they are both zero.
So $$\chi(\omega_X)
=\chi(R^2a_*\omega _X)=\chi(\omega _Y)=q(X)-1.$$

Still looking at the fibration $a:X\longrightarrow Y$, one
sees that $a$ is relatively minimal since $X$ is
minimal. Therefore $K_{X/Y}$ is nef by Theorem 1.4 of \cite{ohno}.
Thus one has $K_X^3\ge (2q(X)-2)K_F^2\ge 2\chi(\omega_X)$, which
is stronger than the required inequality.
\end{proof}

\section{\bf Examples}
In Example 2(e) of \cite{CC99}, one may find a smooth projective
3-fold of general type which is composed with a pencil of surfaces
of $p_g(F)=5$, the biggest value among known examples. Here we
present another example which is composed with curves of genus
$g=5$.

\begin{exmp}
We follow the Example in \S 4 of \cite{CH2004}. We consider
bi-double covers $f_i:C_i\to E_i$ of curves where,
$g(E_i)=0,0,2$. We assume that
$$(d_{i})_* \OO _{C_i}=\OO _{E_i}\oplus
L_i^\vee \oplus P_i^\vee \oplus L_i^\vee \ot P_i^\vee$$ where for
we have $\deg (L_1)=d_1$, $\deg (L_2)=d_2$, $\deg (P_1)=\deg (P_2)=1$ and
$L_3,P_3$ are distinct $2$-torsion elements in $\Pic ^0(E_3)$. In
particular $g(C_i)=2d_i-1$ for $i\in \{1,2\}$ and $f_3$ is \'etale.
It follows that $$\delta : D_1\times
D_2\times D_3\to E_1\times E_2 \times E_3$$ is a $\Z _2^6$ cover.
We denote by $l_i,p_i,l_ip_i$ the elements of $\Z ^2_2$ whose
eigensheaves with eigenvalues $1$ are $L_i^\vee$, $P_i^\vee$ and
$(L_i\otimes P_i)^\vee$. Let $X:= D_1\times D_2\times D_3/G$ where
$G\cong \Z ^4_2$ is the group generated by $$\{ (1,p_2,l_3), (p_1,
l_2, 1), (l_1,1,p_3), (p_1,p_2,p_3) \}.$$ Then one sees that $X$
is Gorenstein and for the induced morphism $f:X\to E_1\times E_2
\times E_3$, one has
$$f_*\OO _X=(\delta _* \OO _{D_1\times D_2\times D_3})^G\cong
\OO _{E_1}\times \OO _{E_2}\times \OO _{E_3}\oplus$$
$$(L_1^\vee \boxtimes L_2^\vee\ot P_2^\vee\boxtimes P_3^\vee )\oplus
(P_1^\vee \boxtimes L_2^\vee\boxtimes L_3^\vee \ot P_3^\vee)
\oplus (L_1^\vee \ot P_1^\vee \boxtimes P_2^\vee\boxtimes
L_3^\vee) .$$ Since $f_* \omega _X = \omega _{E_1\times E_2\times
E_3}\ot f_* \OO _X$, it follows easily that $$H^0(\omega _X)\cong
H^0(\omega _{E_1}\ot L_1)\otimes H^0(\omega _{E_2}\ot L_2\ot
P_2)\otimes H^0(\omega _{E_3}\ot P_3).$$ In particular $p_g(X)=(d_1-1)(d_2-1)$
and $\varphi _1$ factors through the map $X\to C_1/\Z _2\times C_2/\Z_2$.  The
fibers of $\varphi _1$ are then isomorphic to $C_3$ and hence have
genus $5$.
\end{exmp}

\begin{question}\label{Q} A very natural open problem is to find sharp
upper bounds of the invariants of $F$ in Theorem \ref{main}. It is
very interesting to find a new example $X$ which is a Gorenstein
minimal 3-fold of general type such that $\Phi_{|K_X|}$ is of
fiber type and that the generic irreducible component in a general
fiber has larger birational invariants.
\end{question}

We remark that the above question is still open in the surface
case. So Question \ref{Q} is probably quite difficult. A first step
should be to construct new examples with bigger fiber invariants.


\end{document}